\documentclass[12pt, reqno]{amsart}
\usepackage{amsmath, amsthm, amscd, amsfonts, amssymb}
\usepackage[bookmarksnumbered, colorlinks, plainpages, ]{hyperref}
\usepackage{mathtools}
\usepackage{blkarray}
\usepackage{multirow} %mesclar linhas de tabelas \multirow
\usepackage{enumerate} % pacote para enumerar automaticamete os i
\usepackage{amsmath, amsthm, amscd, amsfonts, amssymb}
\usepackage[bookmarksnumbered, colorlinks, plainpages, ]{hyperref}
\usepackage[all]{xy}
\usepackage{pgfplots}
\usepgfplotslibrary{dateplot}

\usepackage{xspace}

\usepackage{color}
\definecolor{darkgreen}{rgb}{0,0.51,0.11}
\usepackage{float}
\usepackage{setspace}
\usepackage{here}
\usepackage{graphicx,float}
\usepackage{psfrag}
\usepackage{mathrsfs}
\usepackage[T1]{fontenc}
\usepackage[utf8]{inputenc}
\usepackage{hyperref}
\usepackage{color}
\definecolor{darkgreen}{rgb}{0,0.51,0.11}
\usepackage{lmodern}
\usepackage{tikz-cd}
\newtheorem{theorem}{Theorem}[section]
\newtheorem{lemma}[theorem]{Lemma}
\newtheorem{proposition}[theorem]{Proposition}

\theoremstyle{definition}
\newtheorem{definition}[theorem]{Definition}
\newtheorem{example}[theorem]{Example}

\theoremstyle{remark}
\newtheorem{remark}[theorem]{Remark}
\numberwithin{equation}{section}

\begin{document}
\setcounter{page}{1}
\title{S-Expansiveness and Zip Shift Maps in Symbolic Dynamics }

\author{S. Lamei }
\address{Faculty of Mathematical Sciences-University of Guilan}
\email{lamei@guilan.ac.ir}

%\urladdr{http://www.icmc.sc.usp.br/$\sim$tahzibi}
\author{ P. Mehdipour}
\address{Departamento de Matematica-Universidade Federal de Viçosa} 
\email{pouya@ufv.br}
\author{ W. Vargas}
\address{Departamento de Matematica-Universidade Federal de Viçosa} 
\email{walterhv@ufv.br}
%\urladdr{http://www.fing.edu.uy/$\sim$ures}

\thanks{}

\keywords{local homeomorphisms, extended symbolic dynamics, shadowing property, expansivity}

\subjclass{Primary: 37Dxx. Secondary: 37D20, 37D45, 26A18, 28D80,   28D05.}

\renewcommand{\subjclassname}{\textup{2010} Mathematics Subject Classification}

\vspace{-0.5cm}
\begin{abstract}
In this paper, we introduce the concept of S-expansiveness for local homeomorphisms and demonstrate that a class of extended symbolic dynamics, known as zip shift maps, are S-expansive and possess the shadowing property. Furthermore, we prove that any S-expansive local homeomorphism is a factor of a zip shift map.

\end{abstract} \maketitle
\vspace{-1cm}
\section{\textbf{Introduction}}

Topological dynamics, which emerged between 1920 and 1930, studies the orbit trajectories  within the framework of their topological behavior. After the 1960, research on the topological dynamics of homeomorphisms \cite{AH} led to significant advances and played a crucial role in the study of differentiable dynamical systems and ergodic theory \cite{K}.

Key properties studied in topological dynamics include the shadowing property and expansivity of homeomorphisms on compact metric spaces. The latter indicates sensitivity to initial conditions, which is a hallmark of chaotic systems. Notably, for expansive systems, the shadowing property is equivalent to topological stability, and all such properties are preserved under topological conjugacy.

One of the earliest examples of expansive homeomorphisms exhibiting the shadowing property is the two-sided shift maps defined on shift spaces \cite{AH}, \cite{MMP}.
The shift maps are suitable examples in which demonstrate the interplay between expansiveness and the shadowing property in dynamical systems. 
The topological conjugacy (or semi-conjugacy) with shift maps serves as a primary tool for investigating these properties within dynamical systems theory. 

In \cite{LM1,LM2}, the authors extend the two-sided shift spaces into a new construct called zip shift space. 
The associated map, referred to as the zip shift, serves as a local homeomorphism and shows significant potential in preserving key properties of topological dynamics. In this paper, we introduce the concept of bilateral S-expansivity and demonstrate that zip shift maps are not only S-expansive but also exhibit the shadowing property. Moreover, we show that S-expansive m-to-1 local homeomorphisms are, in fact, factors of zip shift maps. We hope that these initial properties highlight the robustness and significance of zip shifts in the study of topological dynamics, making them a valuable tool for further exploration and understanding of the dynamical behavior of non-invertible maps.

\section{Preliminary and Main results}

Let $(X,d)$ be a compact metric space and $f:X\to X$ a continuous map. We set, $X^f=\{\tilde{x}=(x_n)|f(x_{n})=x_{n+1},\,n\in\mathbb{Z}\}$ as the inverse limit space corresponded to $X$ and $\tilde{f}:X^f\to X^f$ as its (bilateral) natural extension.  
%There exists a natural projection $\pi: X^f\to X$ such that for any $\tilde{x}=(x_n)\in X^f$, $\pi(\tilde{x})=x_0=x$. 
The future orbit of a point over $f,$ the $O^+(x)=\{f^n(x)|n\in \mathbb{Z}^+\}$ is unique, but with infinitely many, different pre-histories, in which creates infinitely many points in $\pi^{-1}(x)\subset X^f.$ Let $d(.,.)$ represent the (Riemannian) metric on $X$. For $\tilde{x}=(x_n), \tilde{y}=(y_n)\in X^f$, we consider $\tilde{d}$ as a natural metric on $X^f$.

$$\tilde{d}(\tilde{x},\tilde{y}):=\sum_{n=-\infty}^{\infty}2^{-|n|}d(x_{n}, y_{n}).$$
The metric topology induced by $\tilde{d}$ and the product topology are equivalent on $X^f$. %Note that the collection of all sets of the form $\pi^{-1}(U)$ 
\begin{definition}[\textbf{m-to-1 local homeomorphism}]\label{def:2.2}
We say that, the continuous map $f:X\to X$ is an m-to-1 local homeomorphism, whenever $X$ is a compact metric space and there exist $m\in\mathbb{N}$ disjoint connected subsets $X_i\subset X$ that cover $X$ with $m$ \textit{local maps} $f_i:X_i\to X$  such that (for all $i=1,\dots,m$) $f_i$ is a homeomorphism.  
	%Moreover if all $X_i$ has the same volume, then we call $f$ a conservative m-to-1 local diffeomorphism.
\end{definition}

Once working with closed Riemannian manifolds, the m-covering maps are examples of such local homeomorphisms. We have the following Theoem from Eilenberg \cite{AH}.
\begin{theorem}\label{thm:2}
Let $X$ be compact and $f: X\to X$ a continuous surjection. If $f$ is a local homeomorphism, then there exist two positive
numbers $\lambda$ and $\mu$ such that each $D\subset X$ with diameter less than $\lambda$
determines a decomposition of the set $f^{-1}(D)$ with the following properties:
\begin{enumerate}
	\item $f^{-1}(D)=D_1\cup\cdots\cup D_k;$
	\item  $f$ maps each $D_i$ homeomorphically onto $D$,
	\item  if $i\neq j $ then no point of $D_i$ is closer than $2\mu$ to a point of $D_j$;
	\item for every $\theta>0$ there exists some $0 < \epsilon < \lambda$ such that if $diam(D)<\epsilon$ (diameter), then $diam(D_j)<\theta$ for all $i=1,\cdots,k$.
\end{enumerate}
If, in addition, $X$ is connected, then there is a constant $m > 0$ such that $f : X\to X$ is an m-to-1 map.
\end{theorem}

We may use the following well known Lemma in the proof of Theorems. Let $(X,d)$ be a metric space. Then by diameter of a subset $A$ of $X$, we mean: $diam(A)=\sup\{d(a,b):a,b\in A\}.$

\begin{lemma}[\textbf{Lebesgue number}]\label{lem:Leb}\cite{AH}
Let $X$ be a compact metric space and $\alpha$ a finite open cover of
$X$. Then there exists some $\delta=\delta(\alpha)$ such that for any $A\subset X$ if the diameter
of $A$ is less than $\delta$ then $A\subset B$ for some $B\in \alpha$. Such a number is called the
Lebesgue number of $\alpha$.
\end{lemma}

\begin{definition}[\textbf{Principal Domain}]\label{PD}
	We say that $Y\subset X$ is a \textit{principal domain} when $\phi(Y)= X$ and $\phi_{|_Y}:Y\to X$ is a homeomorphism.
	% that $\overline{Y}$ contains only one fixed point of $\phi$.
\end{definition}

\begin{definition}\label{PT}
	Let $(X,d)$ be a compact metric space. We say that a finite cover of $X$ denoted by $\mathcal{D}=\{D_0,D_1,...,D_{N-1}\}$ is a \textit{topological partition} if:
	\begin{enumerate}
		\item Each $D_i$ is open.
		\item $D_i\cap D_j=\emptyset, i\neq j$.
		\item $X = \overline{D_0}\cup\overline{D_1}\cup...\cup\overline{D_{N-1}}$.
	\end{enumerate}
\end{definition}

\begin{definition}[\textbf{Domain topological  partition} \cite{ML}]\label{def:3}
	Consider an m-to-1 local homeomorphism $\phi$ defined on a compact connected metric space $X$ with principal domains $P_i,i=1,\dots,m$. Let $\{P_{i_1},...,P_{i_{k_i}}\}, k_i\in \mathbb{N}$ be a topological partition for $P_i$. We define the \textit{Domain Topological Partition (DTP)} as the collection of all elements in  $\mathcal{P}=\{P_{i_1},...,P_{i_{k_i}}\}_{i=1}^{m}$. Note that by Definition \ref{PD}, for all $k_j\in \mathbb{N}, i=1,2,\dots, m,$ the map $f:P_{i_{k_j}}\to f(P_{i_{k_j}})$ is a homeomorphism.
	%We may refer to this partition as the \textit{Principal Domain Partition or PDP}.
	%We call such a partition \textbf{Main Topological Partition}.
\end{definition}
\begin{definition}[\textbf{Image topological partition} \cite{ML}]\label{def:4}
	Let $\phi: X \rightarrow X$ be an m-to-1 local homeomorphism with a domain topological partition $\mathcal{P}=\{P_{i_1},...,P_{i_{k_i}}\}_{i=1}^{m}$ ($k_i\in \mathbb{N}$). We say that $\mathcal{Q}=\{\mathcal{Q}_1,...,\mathcal{Q}_k\}$ is an \textit{Image topological partition (ITP)} associated to $\mathcal{P}$, if $\mathcal{Q}_j \in \bigvee_{i_j} \phi(P_{i_{j}})$ for $P_{i_j}\in \mathcal{P}$ is a topological partition for $X$ (by definition $\phi(P_i)=X$ when $P_i$ is a principal domain).
\end{definition}

\begin{definition}[\textbf{Extended topological partition} \cite{ML}]\label{def:5}
	Let $\phi: X \rightarrow X$ be an m-to-1 local homeomorphism with an image topological partition $\mathcal{Q}=\{\mathcal{Q}_1,...,\mathcal{Q}_k\}$. The pre-image of $\mathcal{Q}$, induces a domain topological partition which we call \textbf{Extended Topological Partition (ETP)} associated to $\mathcal{Q}$.
\end{definition}

\begin{remark}\label{rem:1}
Note that by definition, m-to-1 local homeomorphisms always admit some ITP induced from some DTP. Indeed it is not difficult to find ETP for such maps.
\end{remark}

The following Definition is adapted from \cite{AH} and \cite{W}.
\begin{definition}[\textbf{Generator/Weak generator}]
Let $f:X\to X$ be a local homeomorphism of a compact metric space. We say that a finite open cover $\alpha$ of $X$ is a \textit{generator} (\textit{weak generator}) for $f$ if 
 for every bi-sequence $\{A_n\}$ of elements of $\alpha,$ the
%and $\tilde{x}=(x_n)\in \pi^ {-1}(A_n)\subset X^f$, that 
$\bigcap_{n=-\infty}^{\infty} f^{-n}%f_{|_{\tilde{x}}}^{-n}
(\overline{A}_n),$ ( $\bigcap_{n=-\infty}^{\infty} f^{-n}(A_n)$) is at most one point. Here $\overline{B}$ denotes the closure of a subset $B$ and as $f$ is non-invertible, $f^{-n}(x)=
(f^{n})^{-1}(x)$. 
%( $\bigcap_{n=-\infty}^{\infty} f_{|_{\tilde{x}}}^{-n}(A_n)$ is at most one point). Here $\overline{B}$ denotes the closure of a subset $B$. 

Equivalently,  one can say that $\alpha$ is a \textit{generator} if for any $\epsilon>0$ 
%and $\tilde{x}=(x_n)\in \pi^ {-1}(A_n)\subset X^f$ 
exists $N>0$ such that the cover $\bigvee_{n=-N}^{N}f^{-n}(\alpha)$;  consists of open sets $A_n$ each of which with diameter at most $\epsilon$, i.e. $\sup_i\,\,diam(A_n)<\epsilon$. %($\bigvee_{n=-N}^{N} f^{-n}(A_n);\, A_n\in \alpha,$ consists of open sets each of which with the diameter at most $\epsilon$). 
	
\end{definition}

\begin{definition}\label{fin P}
	Let $\mathcal{\alpha}=\{A_1,\dots,A_k\}$ and $\mathcal{\beta}=\{B_1,\dots, B_m\}$ be two topological partitions. We say that $\mathcal{\beta}$ is finer than $\mathcal{\alpha}$ (or $\beta$ is a refinement of $\alpha$) if $\#(\mathcal{\alpha})< \#(\mathcal{\beta})$ and for any $B_i\in \mathcal{\beta}$ there exists some $A_j\in \mathcal{\alpha}$ such that $B_i\subseteq A_j$.
\end{definition}

For continuous non-invertible maps, it is commonly used the notions of positive expansivity and c-expansivity \cite{AH}. The latter means for any $\tilde{x}=(x_i), \tilde{y}=(y_i)\in X^f$ if $d(x_i, y_i)<c$ for
all $i\in\mathbb{Z}$ then $\tilde{x}=\tilde{y}$. 
\begin{definition}
	Let $(X,d)$ be a compact metric space and $f:X\to X$ a continuous dynamical system. 
    
    \begin{itemize}
    \item We say that $f$ has \textbf{Positive-Expansivity} if there exists some $\delta>0$ such that for any $x,y\in X$ there exists some $n\in\mathbb{N}$ such that $ d(f^n(x),f^n(y))>\delta$.
    \item We say that $f$ has \textbf{S-Expansivity} and $\gamma>0$ is an \textit{Special Expansivity constant} for $f$ if for any $x,y\in X$ 
 %exists $\tilde{z}\in X^f$  with $\pi(\tilde{x})=x, \pi(\tilde{y})=y$ 
exists some $n\in\mathbb{Z}$ that, $ d(f^n(x),f^n(y))>\gamma$. 
 In case of $n$ being a negative integer $$ d(f^n(x),f^n(y)):=d((f^n)^{-1}(x),(f^n)^{-1}(y)).$$
 Here by $(f^n)^{-1}(x)$ we mean the special branch of $n$ forward local dynamics.
    \end{itemize}
    
\end{definition}

\begin{remark}\label{rem:2}
Note that the c-expansivity is weaker than the positive expansivity and S-expansivity for local homeomorphisms is, in fact, an extension of two-sided expansivity for homeomorphisms. If the map is a homeomorphism, then S-expansivity is equivalent to the usual two-sided expansivity.
\end{remark}

\begin{definition} \label{def:4.7}  
Let $(X,d)$ be a compact metric space and $f:X\to X$ a continuous dynamical system (local hommeomorphism).
\begin{itemize}
    \item We say that the bi-sequence $\{x_n\}_{-\infty}^{+\infty}\subset X$ is a \textbf{$\delta$-pseudo orbit} if $$d(f (x_n),x_{n+1})<\delta.$$ 
    \item We say that $\delta$-pseudo orbit  $\{x_n\}_{-\infty}^{+\infty}$ is \textbf{$\epsilon$-traced} by $y$, if it is $\epsilon$-traced by some $y\in X$, i.e. $d(f^n(y),x_n)<\epsilon$ for all $n\in\mathbb{Z}$.
    \item We say that $f$ has the \textbf{Shadowing property or P.O.T.P.} if for any $\epsilon>0$ exists some $\delta>0$ that any  $\delta$-pseudo orbit of $f$ is $\epsilon$-shadowed by some $y\in X$.
\end{itemize}

\end{definition}

\begin{definition}
	Let $(X,d)$ be a compact connected metric space and $f:X\to X$ a continuous dynamical system. 
    
    \begin{itemize}
    \item We say that $f$ is \textbf{Transitive} if there exists $x\in X$ such that the forward orbit $O^+_{f}(x)$ is dense in $X$.
    \item We say that $f$ is \textbf{Pre-transitive} if there exists $x\in X$ such that the set of backward orbit $O^-_{f}(x)=\{f^{-n}(x):n\geq 0\}$ is dense in $X$.

    \item We say that $f$ is \textbf{Topologically transitive} if for any $U,V$ open subsets of $M$, there exists some $n>0$ such that $f^{n}(U)\cap V \neq \emptyset$.
    \item We say that $f$ is \textbf{Topologically pre-transitive} if for any $U,V$ open subsets of $M$, there exists some $n<0$ such that $f^{n}(U)\cap V \neq \emptyset$.
    
    \end{itemize}
    
\end{definition}

\begin{lemma}\label{pre-trans}
	Let $M$ be a compact metric space (without isolated points) and $f:M\to M$ an onto continuous map. Then Transitivity implies the Topological transitivity and pre-transitivity implies the Topological pre-transitivity.	
\end{lemma}

\begin{proof}
	As the first claim is a well known fact we just prove the secound part \cite{KH}. Let $U, V$ be two disjoint non-empty open subsets of $M$. Then by pre-transitivity, there exists positive integers $k_1,k_2$ such that $f^{-k_1}(x)\in U$ and $f^{-k_2}(x)\in V$. One can without any loss of generality assume that $k_1<k_2$. Then $f^{-k_3}(U)\cap V\neq \emptyset,$ where $k_3=k_2-k_1$.
\end{proof}

\begin{definition}[\textbf{Sensitive dependence on initial conditions}]
	Let $(X,d)$ be a compact metric space and $f:X\to X$ a continuous dynamical system. We say that $f$ has "\textit{sensitive dependence on initial conditions}", if there is $\delta>0$ such that for each $z\in X$ and each
	neighborhood $U$ of $z$ there exist $y\in U$ and $n\in\mathbb{Z}$ such that $d(f^n(z),f^n(y))>\delta$. From definition it follows that X has no isolated points.
\end{definition}
In \cite{AH} the authors show the following for hommeomorphisms. Here we show it for $f$ being a local homomorphism and $X$ a compact metric space.
For $n>0$, let $\textrm{Per}_n=\{x\in X;f^n(x)=x\}$ and $\textrm{Per}(f)=\displaystyle\bigcup_{n=1}^{\infty}\textrm{Per}_n$.
%Rem 2.2.2 Aoki
\begin{proposition}\label{prop:2}
	Let $f:X\to X$ be a local homeomorphism of a compact
	metric space. Suppose $X$ is an infinite set. If $f$ is topologically transitive and  $\textrm{Per}(f)$ is dense in $X$, then $f$ has sensitive dependence on initial conditions.
\end{proposition}
\begin{proof}
%%%%%% walter------->
%\begin{enumerate}
%\item 
We first prove that $\textrm{Per}(f)\neq X$. Otherwise, there exists $n>0$ 
such that $\textrm{int}(\textrm{Per}_n)\neq\emptyset$. Since $f$ is topologically transitive, then $\textrm{Per}_n=X$ and $X$ must be a finite set 
which is a contradiction.

Now we claim that there is $\delta_0>0$ such that for any $x\in X$, there exists
$p\in \textrm{Per}(f)$ such that $d(x,\mathcal{O}_f(p))>\delta_0$. Otherwise, 
%for any $n>0$ there exists a sequence $x_n\in X$ such that 
let  $q_1$, $q_2$ be two arbitrary disjoint periodic points. Let $\delta_0=1/2 (\mathcal{O}_f(q_1),\mathcal{O}_f(q_2))$. For any $x\in X$ with 
$d(x , \mathcal{O}_f(q_1))<\delta_0$, then  $d(x , \mathcal{O}_f(q_2))>\delta_0$ which is a contradiction. 

We will proof that $f$ has sensitive dependence on initial conditions with sensitivity constant $\delta=\delta_0/4$.

Let $x$ be an arbitrary point in $X$ and $N$ be some neighborhood of $x$.
 Since the periodic points are dense, there exists a periodic point $p$ of period $n$ in  $U=N\cap B(x,\delta)$, where  $B(x,\delta)$ is a ball of radius $\delta$. 
As we showed above, there exists a periodic point $q\in X$ satisfying  $d(x, \mathcal{O}_f(q))>4\delta$. 
Let
$$V=\displaystyle\bigcap_{i=1}^{n}f^{-i}(B(f^{i}(q,\delta)))$$

As $f$ is a local homeomorphism, $V$ is open and non-empty. Consequently, since $f$ is transitive, there exists $y\in U$ and a natural number $k$ such that $f^k(y)\in V$.

Now let $j>0$ be the integer part of $\dfrac{k}{n}+1$. So $1\leq nj-k\leq n$. By construction, one has
$$f^{nj}(y)=f^{nj-k}(f^k(y))\in f^{nj-k}(V)\subset B(f^{nj-k}(q),\delta).$$ 
Now $f^{nj}(p)=p$,  $$d(f^{nj}(p),f^{nj}(y))=d(p,f^{nj}(y))$$ and 
$$d(x,f^{nj-k}(q))\leq d(x,p)+d(p,f^{nj}(y))+d(f^{nj}(y),f^{nj-k}(q)).$$
Therefore, since $p\in B(x,\delta)$, 
$$d(f^{nj}(p),f^{nj}(y))>4\delta-\delta-\delta=2\delta.$$
Thus, using triangle inequality again, we have either $d(f^{nj}(x),f^{nj}(y))>\delta$ or  $d(f^{nj}(x),f^{nj}(p))>\delta$. 
%In either case, we have found a point $y\in N$ where $d(f^{nj}(p),f^{nj}(y))>2\delta$. 
This completes the proof.
%\end{enumerate}
\end{proof}

It is easy to verify that topological pre-transitivity implies the topological transitivity. Indeed Proposition \ref{prop:2} is valid if one considers pre-transitivity instead of transitivity.
\begin{proposition}
	Let $f:X\to X$ be a local homeomorphism of a compact
	metric space. Suppose $X$ is an infinite set. If $f$ is topologically pre-transitive and the set $Per(f)$ is dense in $X$, then $f$ has sensitive dependence on initial conditions.
\end{proposition}
\begin{theorem}\label{thm:1}
	The following properties is valid for S-expansive maps.
\begin{enumerate} 
\item If $f: X\to X$ is S-expansive and $Y$ is a closed invariant subset of $X$ (i.e. $f^{\pm 1}(Y)=Y$),
then $f_{|_{Y}}:Y\to Y$ is S-expansive.
\item If $f_i:X_i\to X_i$ ($i = 1,2$) are S-expansive, then the local homeomorphism
$f_1\times f_2: X_1\times X_2\to X_1\times X_2$ defined by
$f_1\times f_2(x,y)=(f_1(x),f_2(y))$
is S-expansive. Every finite direct product of expansive homeomorphisms is expansive.
\item If $X$ is compact and $f:X\to X$ is S-expansive, $h\circ f\circ h^{-1}:Y\to Y$
is S-expansive where $h: X\to Y$ is a homeomorphism.
\end{enumerate}
\end{theorem}

\begin{proof}
The items $(1)$ and $(2)$ have a trivial proof. We prove  item $(3)$. Let $g=h\circ f\circ h^{-1}$. It is not difficult to verify that, $h^{-1}\circ g^n=f^n\circ h^{-1}$ for all $n\in\mathbb{Z}$. Let $\gamma$ be the S-expansiveness constant of $f$. If there exists $x,y \in X$, such that %$\tilde{x}\in \pi^{-1}(x),\tilde{y}\in \pi^{-1}(y)$, 
$d(f^n(x), f^n(y))< \gamma$, for all $n\in\mathbb{Z}$, then $x=y$. By uniform continuity of $h$ and $h^{-1}$, for all $\epsilon> 0,$ in particular $\epsilon=\gamma,$ there exists $\delta>0$ such that $d(g^n(x),g^n(y))<\delta$ implies $d(h^{-1}\circ g^n(x),h^{-1}\circ g^n(y))<\gamma$. 

We have $$d(h^{-1}\circ g^n(x),h^{-1}\circ g^n(y))=d(f^n\circ h^{-1}(x),f^n\circ h^{-1}(y))<\gamma$$
for all $n\in\mathbb{Z}$. Therefore, $h^{-1}(x)=h^{-1}(y)$, which implies $x=y$.
\end{proof}

\subsection{Zip shift spaces}

In this section, we describe the zip shift maps\cite{LM2}. 
Let $Z=\{a_1,a_2,\cdots, a_k\}$ and $S=\{0,1,\dots, N-1\},$ be two collections of symbols that $N\geq k$ and  $\tau:S\to Z$ a surjective map. 
Consider $\prod_{-\infty}^{+\infty}S$ and let $\bar t=(t_i)\in\prod_{-\infty}^{+\infty}S$. Then to any such point $\bar t=(t_i)_{i\in \mathbb{Z}}$ correspond a point $\bar x=(x_i)_{i\in \mathbb{Z}}$, such that
\begin{equation}\label{Z}
x_i=\left\{\begin{tabular}{ll}
$t_i\in S \hspace{7mm} $\,\,\,\,\quad\quad\quad$\forall i\geq 0$\\
$\tau(t_{i})\in Z \hspace{7mm} $\,\,\,\quad\quad $\forall i<0$. 
\end{tabular}\right.
\end{equation}
Define 
$\Sigma= \{\bar x=(x_i)_{i\in \mathbb{Z}}\ : x_i \textrm{\ satisfies} \,(\ref{Z})\}.$
We define $\sigma_{\tau}:\Sigma\to\Sigma$ as follows.
\begin{eqnarray}\label{ZS}
(\sigma_{\tau}(\bar x))_i=\left\{\begin{tabular}{ll}
$x_{i+1} \,\,\,\quad\quad \text{if}\,\,i\neq -1,$ \\
$\tau(x_{0}) \quad\quad \text{if}\,\,i=-1.$
\end{tabular}\right.
\end{eqnarray}
We call this map \textit{zip shift}. If $\Sigma$ is closed and invariant under $\sigma_{\tau}$, call it \textit{Zip shift Space}.
One can equip $\Sigma$ with a metric $d$. 
%For  $\bar{s},\,\bar{t}\in \Sigma$, let define \begin{eqnarray}
%N^+=N^+(\bar{s},\,\bar{t})=1+\min_{i} \{ i\geq0| s_i\neq t_i, i\in \mathbb{Z}\},\\
%N^-=N^-(\bar{s},\,\bar{t})=\min_{i} \{ i\leq0| s_i\neq t_i, i\in \mathbb{Z}\}.
%\end{eqnarray}
%Then$\bar d(s,\,t)=\frac{1}{2}(\frac{1}{2^{N^+}}+\frac{1}{2^{N^-}})$ defines a metric on $\Sigma$ 
Let $M:\Sigma\times \Sigma\rightarrow\mathbb{N}\cup \{0\}$ be given as, 
\begin{align}\label{metric}
    M(x,y)=\left\lbrace\begin{array}{lr}
\infty,&\text{ if } x=y\\
\min\{|i|;\,\{x_i\}\neq \{y_i\}\},&\text{ if }x\neq y
\end{array}\right.
\end{align}
then  $d(x,y)=\frac{1}{\lambda^{M(x,y)}}$ defines a metric on $\Sigma$, where we choose $\lambda=\#(S)$. 
Throughout the remainder of this paper, we work within the metric space$(\Sigma, d)$.

\begin{proposition}\cite{MM}
	The followings are valid for $\sigma_\tau:\Sigma \to \Sigma$.
	\begin{itemize}
		\item $\sigma_{\tau}(\Sigma)=\Sigma$ (full invariant).
		\item $\sigma_{\tau}$ is a local homeomorphism.
	\end{itemize}
\end{proposition}

\begin{proposition}\cite{LM1}
	Periodic points of a zip shift map are dense in the zip shift space. $Per(\sigma_{\tau}^k)=(l)^{k}$. Moreover the zip shift map is topologically mixing.
\end{proposition}

The basic cylinders on a zip shift space are defined as follows.
\begin{equation}\label{B-C}
C_i^{s_i}=\{x\in\Sigma\,|\,x_i=s_i, i\in\mathbb{Z},\,s_i\in Z,\, \text{if}\,i<0,\,\text{and}\,s_i\in S,\, \text{if}\,i\geq 0\}.
\end{equation}
The set $C_i^{s_i}$ presents the set of all sequences, that have $s_i$ in the $i-$th entry.  For $i,n\in \mathbb{Z}$ and $\ell\in \mathbb{N}\cup \{0\}$, one can define a general cylinder set as follows. 
$$C_{i_1\dots,i_k}^{s_{1},\dots,s_{k}}=\{ (t_n)\in \Sigma|t_{i_1}=s_1,\cdots,t_{i_k}=s_k; s_1,\dots,s_k\in S\cup Z\},$$
where  $s_j\in Z,\, \text{if}\,i_j<0,\,\text{and}\,s_j\in S,\, \text{if}\,i_j\geq 0$ ($1\leq j\leq k$). 
  
Such cylinder sets are clopen subsets in the metric topology of $(\Sigma, d)$.  The set of all such cylinder sets, form a basis for the topology. In fact, $(\Sigma,d)$ is compact, totally disconnected and perfect, which means it is a Cantor set.

\begin{example}
	Let $S=\{0,1,2,3\}$ and $Z=\{a,b\}$. Define the corresponding transition map $\tau:S\to Z$ as $\tau(0)=\tau(2)=a$ and $\tau(1)=\tau(3)=b$.  
	Let $\bar x=(x_n)_{n\in\mathbb{Z}}=(\cdots a\,b\,a\,b\,b\,\textbf{.}\,1\,0\,3\,1\,1\,2\,\cdots).$  One can verify that $$\sigma_{\tau}((\cdots\,a\,b\,a\,b\,b\,\textbf{.}\,1\,0\,3\,1\,1\,2\,\cdots))=(\cdots\,a\,b\,a\,b\,b\,b\,\textbf{.}\,0\,3\,1\,1\,2\,\cdots).$$

 Note that one can consider $S=S_1\cup S_2$, where, $S_1=\{0,1\}, S_2=\{2,3\}$ and define $\tau$ as follows.
 \begin{eqnarray}\label{TL}
\tau(x)=\left\{\begin{tabular}{ll}
$\tau_1(x) \,\,\,\quad\quad \text{if}\,\,x\in S_1,$ \\
$\tau_2(x) \quad\quad \text{if}\,\,x\in S_2.$
\end{tabular}\right.
\end{eqnarray}
Then, $\sigma_{\tau_i}:\Sigma_i\to \Sigma$ will be the local dynamics of $\sigma_{\tau}$ (see Definition \ref{def:2.2}) with $\Sigma_i=C_{0}^{s_m}\cup C_0^{s_n},s_m\neq s_n\in S_i$  for $i=1,2$.$\hfill\blacksquare$
		
\end{example} 

\begin{example}\label{Ex:2}
	Let $X=[0,1]$ and $T(x)=3x\,\,mod\,1$ be the expanding map. For $S=\{0,1,2\},\, Z=\{a\}$, this map is conjugate (mod 0) with a zip shift map with $\tau:S\to Z$ that $\tau(0)=\tau(1)=\tau(2)=a$. $\hfill\blacksquare$
%	
%	\begin{figure}[h]
%		\includegraphics[width=0.35\textwidth]{2x.pdf}
%		\caption{$f(x)=2x\,\, mod 1$}
%		\label{Fig:2}  
%	\end{figure}	 
	
\end{example} 
\begin{proposition}
	The zip shift maps are S-expansive local homeomorphisms.
\end{proposition}

\begin{proof}

Let $\gamma=1/2$. For any $x,y \in \Sigma$ recall that  $d(x,y)=\frac{1}{\lambda^{M(x,y)}}$ as in \eqref{metric}. 
For $x\neq y $ there exists some $i\in \mathbb{Z}$ such that $x_{i}\neq y_{i}$. Let $i$ be the least in modulus of such $i$. If $i>0$, then $n=i+1$ and  $d(\sigma_{\tau}^{i+1}x,\sigma_{\tau}^{i+1}y)=\frac{1}{\lambda^0}>1/2$. If $i<0$, then $d(\sigma_{\tau}^{-i}x,\sigma_{\tau}^{-i}y)=d((\sigma_{\tau}^{i})^{-1}x,(\sigma_{\tau}^{i})^{-1}y)=\frac{1}{\lambda^0}>1/2$. 
Indeed the zip shift map is an S-expansive dynamical system with $\gamma=1/2$.
\end{proof}
\begin{theorem}
	Let $\sigma_{\tau}:\Sigma_{Z,S}\to \Sigma_{Z,S}$ be a full zip shift map. Then $\sigma_{\tau}$ has shadowing property.
\end{theorem}

\begin{proof}
 Let $(\Sigma_{Z,S},\sigma_{\tau})$ be a zip shift space.  Let $\epsilon>0$ and $m\geq 1$ such that $\lambda^{-m}<\epsilon$. Then if $x,y\in \Sigma_{Z,S}$ such that $d(x,y)<\lambda^{-(m+1)}$, we have $x_i=y_i$ for $|i|\leq m$. In particular, $x_0=y_0$. Let $\{z^j\}_{j=-\infty}^{+\infty}$ be a $\delta=\lambda^{-(m+1)}$-pseudo orbit for $\sigma_{\tau}$. Therefore by definition, $d(\sigma_{\tau}(z^j),z^{j+1})<\delta$ and indeed, $ z^j_{k+1}=z^{j+1}_k$ for all $j\in\mathbb{Z}, |k|\leq m$. Let $x=(x_n)$ where 
 $$x=(\dots \,z_{-m-2}^{-2}\,z_{-m-1}^{-1}\,z_{-m-1}^{0}\,\dots\,z_{-1}^{0}\,;z_0^0\,\dots\,z_m^0\,z_{m+1}^{0}z_{m+1}^{1}\,z_{m+1}^{2}\,\dots).$$
 As $\{z^j\}_{j=-\infty}^{+\infty}$ is a $\delta$ pseudo orbit, so for all $|j|\leq m$ and $n\in\mathbb{Z}$, we have $x_{n+j}=z_j^n$, which means $$d(\sigma_{\tau}^n(x),z^n)<\lambda^{-m} <\epsilon.$$ 
 Thus $\sigma_{\tau}$ has the shadowing property.
\end{proof}

In \cite{M} the authors show that  a zip shift map has shadowing property iff it is a zip shift of finite type (SFT-zip shift).

%\section{Main Results}
\begin{theorem}\label{thm:3}
	Let $f:X\to X$ be an m-to-1 local homeomorphism of a compact metric
	space. Then the following are equivalent.
	\begin{enumerate}
		\item $f$ is S-expansive;
		\item $f$ has a generator;
	\item $f$ has a weak generator;
%            \item $f$ is S-expansive;
	\end{enumerate}
	\end{theorem}
	
\begin{proof}
$(1)\Rightarrow (2)$: 
Let $f$ be S-expansive with expansivity constant $\gamma$. Let $\alpha=\{A_1,\dots,A_k\}$ be a finite cover consisting of open balls of radius less than $\gamma/2$. Suppose that $x,y \in\bigcap_{-\infty}^{\infty} f^{-n}(\overline{A}_n)\}$.
As by S-expansivity, $d(f^n(x),f^n(y))\leq\gamma/2<\gamma$, for all $n\in \mathbb{Z}$, therefore   $x=y$. 

$(2)\Rightarrow (3)$: Suppose that $f$ has a generator $\alpha$ which is not a weak generator. Therefore %there exists $\tilde{z}=(z_n)\in \pi^ {-1}(A_n)$ and 
let exists $x,y\in\bigcap_{-\infty}^{\infty} f^{-n}(A_n)\neq \emptyset$ with $x\neq y$. As for any $n,$ the $A_n\subset \overline{A}_n,$ it means that $x,y \in \bigcap_{-\infty}^{\infty}f^{-n}(\overline{A}_n)\neq \emptyset,$ with $x\neq y$, which is a contradiction.

$(3)\Rightarrow (1)$: Let $\beta$ be a weak generator for $f$ with Lebesgue number $\delta$. Therefore, for every bi-sequence $\{B_n\}$ of elements of $\beta,$ the
%and $\tilde{z}=(z_n)\in \pi^ {-1}(B_n)\subset X^f$, the 
$\bigcap_{n=-\infty}^{\infty} f^{-n}(B_n)$ is at most one point. Suppose that $d(f^n(x),f^n(y))<\delta$ for all $n\in \mathbb{Z}$. Then by Lemma \ref{lem:Leb} for each $n$ there is some $B_n\in \beta$ such that
$f^n(x),f^n(y)\in B_n$, and so $x,y\in \bigcap_{-\infty}^{\infty} \{f^{-n}(B_n)\}$ which is at most one point. Thereupon, $x=y$ and $f$ is S-expansive.
\end{proof}

\begin{theorem}
Let $f: X\to X$ be an m-to-1 local homeomorphism of a compact metric
space and let $k\neq 0$ be a positive integer. Then $f$ is S-expansive if and only if $f^k$ is S-expansive.
\end{theorem}

\begin{proof}
Let $\gamma$ be the S-expansivity constant of $f$. Then by Theorem \ref{thm:3} $f$ has a generator. Call it $\alpha$. Note that $\bigvee_{i=0}^{k-1}f^{i}(\alpha)$ will be a generator for $f^k$. Moreover, if  $\alpha$ is a generator for $f^k$, then it is also a generator for $f$. 
\end{proof}

\begin{theorem}[\textbf{Main}]
	Let $f: X\to X$ be an S-expansive m-to-1 local homeomorphism of a compact connected metric space. Then there exist finite symbolic sets $S, Z,$ a closed invariant subset $\Sigma$, a zip shift map $\sigma_{\tau}:\Sigma\to \Sigma$ and a continuous surjection
	$\pi:\Sigma\to X$ satisfying $\pi\circ \sigma_\tau= f\circ \pi$, i.e. the following diagram commutes.
	\begin{equation}\label{1-0}
	\begin{tikzcd}
	[row sep=tiny,column sep=tiny]& & & \\ \Sigma \arrow{rr}{\sigma_{\tau}}\arrow{dd}{\pi} & & \Sigma\arrow{dd}{\pi} \\& \circlearrowright &  & \\ X\arrow{rr}{f} & & X\\
\end{tikzcd}
\end{equation}
Moreover, if $X$ is connected, the above result is valid for any S-expansive local homeomorphism $f: X\to X$.
\end{theorem}
\begin{proof}
Let $\gamma>0$ be the S-expansive constant for $f$. As $f$ is an m-to-1 local homeomorphism, there exist $m$ principal domains $D_i$ such that $f_i=f_{_{|_{D_{i}}}}:D_i\to X$ is a homeomorphism. Note that the collection of principal domains $\mathcal{D}=\{D_1,\dots,D_m\}$ is a topological partition for $X$. We choose a finite topological partition $\alpha=\{A_0,\cdots, A_{k-1}\},$
such that for all $A\in \alpha$, the  $diam(A_i) < \gamma, i=0,\dots, k-1$. Refine this topological partition to $\beta= \alpha \vee D$ if necessary to obtain a DTP. Let $\beta=\{B_0,\cdots,B_{l-1}\}$ and  $\partial_{\infty}=\bigcup\limits_{n=-\infty}^{+\infty} f^{n}(\partial{B_i})$. Then $\partial_{\infty}$ 
is a set of first category and
so $\bar X=X\setminus \partial_{\infty}$ is dense in $X$. Suppose that $\Sigma_{Z,S}$ is a zip shift space in which $S,Z$ and the transition map $\tau$ are deduced by the help of $\beta$ and some induced image topological partition (Remark \ref{rem:1}). We show that for any $x\in \bar X$,  a unique member of $\Sigma_{Z,S}$ is assigned.  Let $Q=\{Q_1,\cdots,Q_s\}$ be the image topological partition associated to $\beta$. Consider the main topological partition $M$, associated with $Q$. If necessary one can refine the partition $Q$ (take $Q\vee \beta$) and $M,$ to obtain a good image and consequently main topological partition with $diam(M_i)<\gamma$ for all $ M_i\in M$. Note that the elements of $M$ are pre-images of $Q$ and by Theorem \ref{thm:3} $M$ is a generator. The good image topological partitions $Q$ and the main topological partitions $M,$ induce two finite symbolic sets, which can be represented by $S$ (for $M$) and $Z$ (for $Q$). Furthermore, as the elements of $M$ are pre-images of the elements of $Q$, there are surjective maps $\hat{\tau}: M\to Q$ and $\tau: S\to Z$ that relate the domain and image partitions, as well as their associated symbolic sets. Let $ \Sigma_{Z,S}$ be the full zip shift map associated with $S,Z$ and with transition map $\tau$. Let $(\#(Z)=\bar s)\leq (\bar k=\#(S))$ and consider $Z=\{a_1,\dots,a_{\bar s} \}$ and $S=\{b_0,b_1,\dots, b_{{\bar k}-1}\}$. Then to any $x\in \bar X$ can be assigned a unique element of the form $(x_n)\in \Sigma_{Z,S}$. Let $M_i\in M$ and $Q_i\in Q$. Then for $n\geq 0,$ $x_n=b_i$ if $x\in f^{-n}(M_i)$ i.e. $f^n(x)\in M_{i}$ and for $n<0,$ $x_n=a_i$ if $x\in f^{-n-1}(Q_{i})$ ($-n>0$). 
%Let $\sigma_{\tau}: \Sigma_{Z,S}\to \Sigma_{Z,S}$ denote the zip shift map with $\tau:S\to Z$ induced from the relation between $\beta$ and $Q$. 
Suppose that $\Lambda\subset \Sigma_{Z,S}$ denotes the set of all such points. In other words, one can say that $x\in \overline{[\bigcap_{-\infty}^{0}f^{-n}(Q_{x_{n-1} })]}\cap\overline{[\bigcap_0^\infty f^{-n}(M_{x_n})]}$. Note that as domain topological partition $M$ is a generator and we have the S-expansivity, the non-empty intersection $\overline{[\bigcap_{-\infty}^{0}f^{-n}(Q_{x_{n-1} })]}\cap\overline{[\bigcap_0^\infty f^{-n}(M_{x_n})]}$ is a single point. One needs to guarantee that $\Lambda$ is a zip shift space. As $\Sigma_{Z,S}$ is a compact metric space and 
$\overline{[\bigcap_{-n+1}^{0}\phi^{-k}(Q_{x_{k-1} })]}\cap\overline{[\bigcap_0^n\phi^{-k}(M_{x_k})]}$ has the finite intersection property, the subset $\Lambda$ of $\Sigma_{Z,S}$ is closed. Let $\sigma_{\tau}: \Lambda\to \Lambda$ and consider the injective map $\phi:\Lambda\to \bar X$, in which $\phi^{-1}(x)=(z_n)_{n\in\mathbb{Z}}$. Note that, $\sigma_{\tau}\circ\phi^{-1}(x)=\phi^{-1}\circ f(x)$, for all $x\in\bar X$. This means that $\Lambda$ is $\sigma_{\tau}-$invariant as well. Indeed, it is a (sub)zip shift space.

In order to extend $\phi$ to a continuous map $\tilde{\phi}:\Lambda\to X $ satisfying $\tilde{\phi}\circ\sigma_{\tau}(z_n)=f\circ\tilde{\phi}(z_n)$ for all $z\in\Lambda$, it is enough to show that $\phi$ is a uniformly continuous map. Let $\epsilon>0$ be arbitrary. As $f$ is an S-expansive m-to-1 local homeomorphism and $M$ is a generator, there exists some $N>0$ such that each member of $\bigvee_{k=-N}^{N}f^{k}(M)$ has diameter less than $\epsilon.$ Moreover, if $a=(a_n),b=(b_n)$ be elements of $\Lambda$ with $a_n=b_n$ for $|n|\leq N$, i.e. for $0<\delta \leq\frac{1}{2^N}$, then $\phi(a)$ and $\phi(b)$ belong to the same element $\overline{\vee_{k=-N+1}^{0}f^{-k}(Q_{x_{k-1}}}) \vee (\overline{\vee_{k=0}^N f^{-k}(M_{x_k})}\in\bigvee_{k=-N}^{N}f^{k}(M)$ and thus $d(\phi(a),\phi(b))<\epsilon.$ Thereupon, $\phi$ is uniformly continuous.
Moreover, the last assertion is correct by Theorem \ref{thm:2}.

\end{proof}

\subsection{Dataavailability}
Data sharing not applicable to this article as no data sets were generated or analysed during the current study.
\subsection{Conflicts of interest}
The authors declare no conflict of interest.
\subsection{Founding}
The second author gratefully acknowledges partial financial support from FAPEMIG, Brazil.

\end{document}